\documentclass[11pt]{amsart} 
\makeatletter
\renewcommand{\@seccntformat}[1]{%
  \ifcsname prefix@#1\endcsname
    \csname prefix@#1\endcsname
  \else
    \csname the#1\endcsname\quad
  \fi}
\newcommand\prefix@section{\thesection: }
\makeatother

\usepackage{amsbsy,amsmath,amsthm,amssymb,color,verbatim,ulem,color,wrapfig}
\usepackage[backrefs]{amsrefs}
\usepackage{float}
\usepackage{bm}
\usepackage{fullpage,cancel}
\usepackage{hyperref}
\usepackage{float}
\usepackage{wrapfig}
\usepackage[dvipsnames]{xcolor}
\usepackage{graphicx}
\usepackage{subcaption}
\usepackage{multicol}
\usepackage{multirow}
\usepackage{longtable,array}
\usepackage{epigraph}
\usepackage{enumitem,titlesec}
\usepackage[mathlines]{lineno}

\usepackage{makecell}
\setcellgapes{4pt}
\usepackage{supertabular}
\usepackage{tikz}
\usepackage{tkz-graph}
\tikzstyle{vertex}=[circle, draw, inner sep=0pt, minimum size=4pt]

\tikzstyle{vtx}=[circle, draw, inner sep=0pt, minimum size=8pt]
\usetikzlibrary{decorations.pathreplacing}

\definecolor{darkgreen}{cmyk}{.9,0,.9,.2}
\definecolor{midgray}{gray}{0.60}
\definecolor{lightgray}{gray}{0.90}
\definecolor{lmgray}{gray}{0.70}
\usepackage{colortbl}
\usetikzlibrary{arrows.meta,calc,decorations.markings,math,arrows.meta}
\usepackage{caption}
\usepackage{subcaption}
\makeatletter
\def\@tocline#1#2#3#4#5#6#7{\relax
  \ifnum #1>\c@tocdepth 
  \else
    \par \addpenalty\@secpenalty\addvspace{#2}%
    \begingroup \hyphenpenalty\@M
    \@ifempty{#4}{%
      \@tempdima\csname r@tocindent\number#1\endcsname\relax
    }{%
      \@tempdima#4\relax
    }%
    \parindent\z@ \leftskip#3\relax \advance\leftskip\@tempdima\relax
    \rightskip\@pnumwidth plus4em \parfillskip-\@pnumwidth
    #5\leavevmode\hskip-\@tempdima
      \ifcase #1
       \or\or \hskip 1em \or \hskip 2em \else \hskip 3em \fi%
      #6\nobreak\relax
    \hfill\hbox to\@pnumwidth{\@tocpagenum{#7}}\par
    \nobreak
    \endgroup
  \fi}
\makeatother

\usepackage{calligra}
\usepackage[T1]{fontenc}
\usepackage{titlesec}

\renewcommand{\thesection}{Part~\arabic{section}}

\usepackage{relsize}
\makeatletter
\newtheorem*{rep@theorem}{\rep@title}
\newcommand{\newreptheorem}[2]{%
\newenvironment{rep#1}[1]{%
 \def\rep@title{#2 \ref{##1}}%
 \begin{rep@theorem}}%
 {\end{rep@theorem}}}
\makeatother

\makeatletter
\newtheorem*{rep@conjecture}{\rep@title}
\newcommand{\newrepconjecture}[2]{%
\newenvironment{rep#1}[1]{%
 \def\rep@title{#2 \ref{##1}}%
 \begin{rep@conjecture}}%
 {\end{rep@conjecture}}}
\makeatother

\makeatletter
\newcommand{\addresseshere}{%
  \enddoc@text\let\enddoc@text\relax
}
\makeatother

\theoremstyle{definition}

\newreptheorem{theorem}{Theorem}
\newrepconjecture{conjecture}{Conjecture}

\newcommand{\N}{\mathbb{N}}

\newcommand{\bref}[1]{\textcolor{blue}{\ref{#1}}}
 
\usepackage{etoolbox}

\allowdisplaybreaks

\title{Parking Functions: Choose Your Own Adventure}
   
\author{Joshua Carlson}
\address[J. Carlson]{Department of Mathematics and Statistics, Williams College}
\email{\textcolor{blue}{\href{mailto:jc31@williams.edu}{jc31@williams.edu}}}

\author{Alex Christensen}
\address[A. Christensen]{University of Arizona, Department of Mathematics}
\email{\textcolor{blue}{\href{mailto:ajc333@comcast.net}{ajc333@comcast.net}}}

\author{Pamela E. Harris}
\address[P. E. Harris]{Department of Mathematics and Statistics, Williams College}
\email{\textcolor{blue}{\href{mailto:peh2@williams.edu}{peh2@williams.edu}}}

\author{Zakiya Jones}
\address[Z. Jones]{Pomona College, Department of Mathematics}
\email{\textcolor{blue}{\href{mailto:zakiyacmjones@gmail.com}{zakiyacmjones@gmail.com}}}

\author{Andr\'es Ramos Rodr\'iguez}
\address[A. Ramos Rodr\'iguez]{University of Puerto Rico, Rio Piedras, Department of Computer Science}
\email{\textcolor{blue}{\href{mailto:ramosandres443@gmail.com }{ramosandres443@gmail.com}}}

\begin{document}

\renewcommand{\abstractname}{Warning}

\begin{abstract}
The reading of this paper will send you down many winding roads toward new and exciting  research topics enumerating generalized parking functions. Buckle up! 
\end{abstract}

\maketitle
\section*{Preface}

Let $n\in\N:=\{1,2,3,\ldots\}$ and consider a parking lot consisting of $n$ consecutive parking spots along a one-way street. Suppose $n$ cars want to park one at a time in the parking lot and each car has a preferred parking spot. Each car coming into the lot initially tries to park in its preferred spot. However, if a car's preferred spot is already occupied, then it will park in the next available spot. Since the parking lot is along a one-way street, it is not guaranteed that every car will be able to park before driving past the parking lot. This dilemma leads to the idea of a parking function.

Let us make this definition precise. For \mbox{$n\in\N$}, let $[n]:=\{1,\dots, n\}$. Formally, suppose the parking spots are labeled $1, 2, \ldots, n$, in order, along the one-way street and the cars are labeled according to the order in which they try to park. In other words, for each $i\in[n]$, car $c_i$ is the $i^{\text{th}}$ car to try to park and prefers spot $a_i\in[n]$. Note that more than one car can have the same preference.  To park, cars first drive to their preferred spot and park in it if it is available. If their preferred spot is occupied then they drive forward and park in the next available spot. If all $n$ cars can park in the parking lot under these conditions, then the preference list $(a_1,a_2, \ldots, a_n)$ is called a {\textbf{parking function}} (of length $n$). For example, $(1,2,4,2,2)$ is a parking function, but $(1,2,2,5,5)$ is not. 
Naturally, the first question that arises is: \textit{``For any $n\in\N$, how many parking functions are there?''} Konheim and Weiss showed that the number of parking functions of length $n$ is $(n+1)^{n-1}$~\cite{KonheimAndWeiss}.

In this paper, enumerating the set of preference lists that allow all cars to park satisfying certain conditions is what we call the \textbf{parking problem.} In the adventure
that follows, we will make a variety of choices to define new parking problems. For instance, we could change the structure of the parking lot, the way cars prefer parking spots, or the way cars drive down the street. Regardless of the choices we make, the quest in the first part of our adventure is to count the number of preference lists that satisfy the conditions of the parking problem at hand. With this enumeration, we then begin a new part of the adventure where we study properties of the set of parking functions arising in these parking problems. This gives rise to an entire new avenue of open problems.
Let us begin!\footnote{This adventure will be an ad free experience. We withhold references until our concluding remarks in \ref{sec:closingremarks}, where we direct the reader to literature related to the problems presented.}

\section{The adventure begins}\label{firstQuestion}
\epigraph{``\textit{Two roads diverged in a wood, and I ---\\ I took the one less traveled by, \\And that has made all the difference.}''}{-{Robert Frost}}

You find yourself alone in your car and at a stand still in a traffic jam with nothing but Stanley's Enumerative Combinatorics Volume 2 to keep you occupied. Having randomly flipped through the book you land on page 94 and find the following problem:

\begin{center}
\begin{minipage}[b]{.85\textwidth}
\cite[\bf 5.29. a.]{EC2} There are $n$ parking spots $1,2,,\ldots,n$ (in that order) on a one-way street. Cars $c_1,c_2,\ldots c_n$ enter in that order and try to park. Each car $c_i$ has a preferred spot $a_i$. A car will drive to its preferred spot and try to park there. If the spot is already occupied, the car will park in the next available spot. If the car must leave the street without parking, then the process fails. If $\alpha=(a_1,a_2,\ldots,a_n)$ is a sequence of preferences that allows every car to park we call $\alpha$ a \textit{parking function}. Show that a sequence $(a_1,\ldots,a_n)\in[n]^n$ is a parking function if and only if the increasing rearrangement $b_1\leq b_2\leq\cdots\leq b_n$ of $a_1,a_2,\ldots,a_n$ satisfies $b_i\leq i$.
\end{minipage}
\end{center}
The traffic is starting to get you down, and your lack of paper and pen makes it impossible to try to solve Stanley's problem. So instead, being the capricious mathematician you are, you decide to ponder possible variations of the parking function problem to pass the time.

First, you need to decide what direction you want the cars to be able to move.
You have three options. If their preferred spots are occupied, do you want the cars to...
 \begin{itemize}
     \item[a.] Move exclusively forward? Go to \bref{onlyforward}.
     \item[b.] Move backward first, then forward? Go to \bref{backwardthenforward}.
     \item[c.] Have variable movement? 
    Go to \bref{variablemovement}.
 \end{itemize}

\subsection{--\ \ Cars only move forward}\label{onlyforward}
$ $\newline

Now, choose something special to further customize your problem.

 \begin{itemize}
     \item[a.] I want the cars in my lot to have variable sizes. Go to \bref{carsizes}.
     \item[b.] I want my parking lot to have some spots that allow multiple cars to park in them. Go to \bref{onespacemanycars}.
 \end{itemize}

\subsection{--\ \ Cars move backward, then forward}\label{backwardthenforward}
$ $\newline

You choose to allow the cars to move backward before moving forward. How many spots back do you want the cars to be able to check?
 \begin{itemize}
     \item[a.] One spot. Go to \bref{onespaceback}.
     \item[b.] Multiple spots. Go to \bref{multiplespacesback}.
 \end{itemize}

\subsection{--\ \ Cars have variable movement}\label{variablemovement}
$ $\newline

Cool, you want to mix things up. Now, only certain cars are allowed to move backward. How do you want this to be decided?
 \begin{itemize}
     \item[a.] I want to specify which cars can move backward. Go to \bref{fixeddirection}.
     \item[b.] I want to make it random. Go to \bref{randomdirection}.
 \end{itemize}

 \subsection{--\ \ Cars move back one spot before moving forward}\label{onespaceback}
 $ $\newline
 
 Each car can now check the spot directly behind it if its preferred spot is occupied. Now, to make it more interesting, do you want the parking lot to have... 
 \begin{itemize}
     \item[a.] Parking spots that are unavailable due to obstructions? Go to \bref{backwardtrailer}.
     \item[b.] Some spots in which multiple cars can park? Go to \bref{backwardmotorcyle}.
 \end{itemize}

\subsection{--\ \ Cars move back multiple spots before moving forward}\label{multiplespacesback}
$ $\newline

Now, you need to decide how the cars move backward. You have two options. 
 \begin{itemize}
     \item[a.] I want the cars to teleport $k$ spots back. Go to \bref{teleportation}.
     \item[b.] I want the number of steps that car $c_i$ can take back to be dependent on the value of $i$. Go to \bref{i-stepsback}.
 \end{itemize}

\subsection{--\ \ Cars of different lengths}\label{carsizes}
$ $\newline

    Suppose you are in the future and each car $c_i$, for $1 \leq i \leq n$, can  teleport and its size $s_i \in \N$ can vary. Sadly, parking technology has not caught up to the car technology and parking spots are only of size $1$. If car $c_i$ is of size $1 \leq s_i \leq n$, it needs $s_i$ consecutive empty parking spots to park. If the car $c_i$ prefers the spot $a_i$, it can park if the spots from $a_i$ to $a_i + s_i - 1$ are empty. If any of these spots are occupied, the car is unable to park and will teleport $k\in\N$ units forward. The way you determine the value of $k$ is up to you. If $a_i + k + s_i - 1 > n$, then the car teleports out of the parking lot and does not park. If $a_i + k + s_i -1 \leq n$, then the car teleports and is still in the lot. The futuristic teleporting car $c_i$ can park if it can find $s_i$ consecutive empty parking spots starting at and including spot $a_i +k$. Unfortunately, each futuristic car has only enough fuel to teleport once, so if the car is unable to park after teleporting, it does not park at all. Can you count the number of preference lists that are parking functions with these futuristic cars?

Finished? Good news! This isn't the end of the road! Once you have counted these parking functions, there are many more properties that you can study. Go to \bref{hungryformore}.

\subsection{--\ \ Multiple cars park in a single spot}\label{onespacemanycars}
$ $\newline
\begin{wrapfigure}{l}{2in}
\includegraphics[trim=0 0 0 0, clip,width=2in]{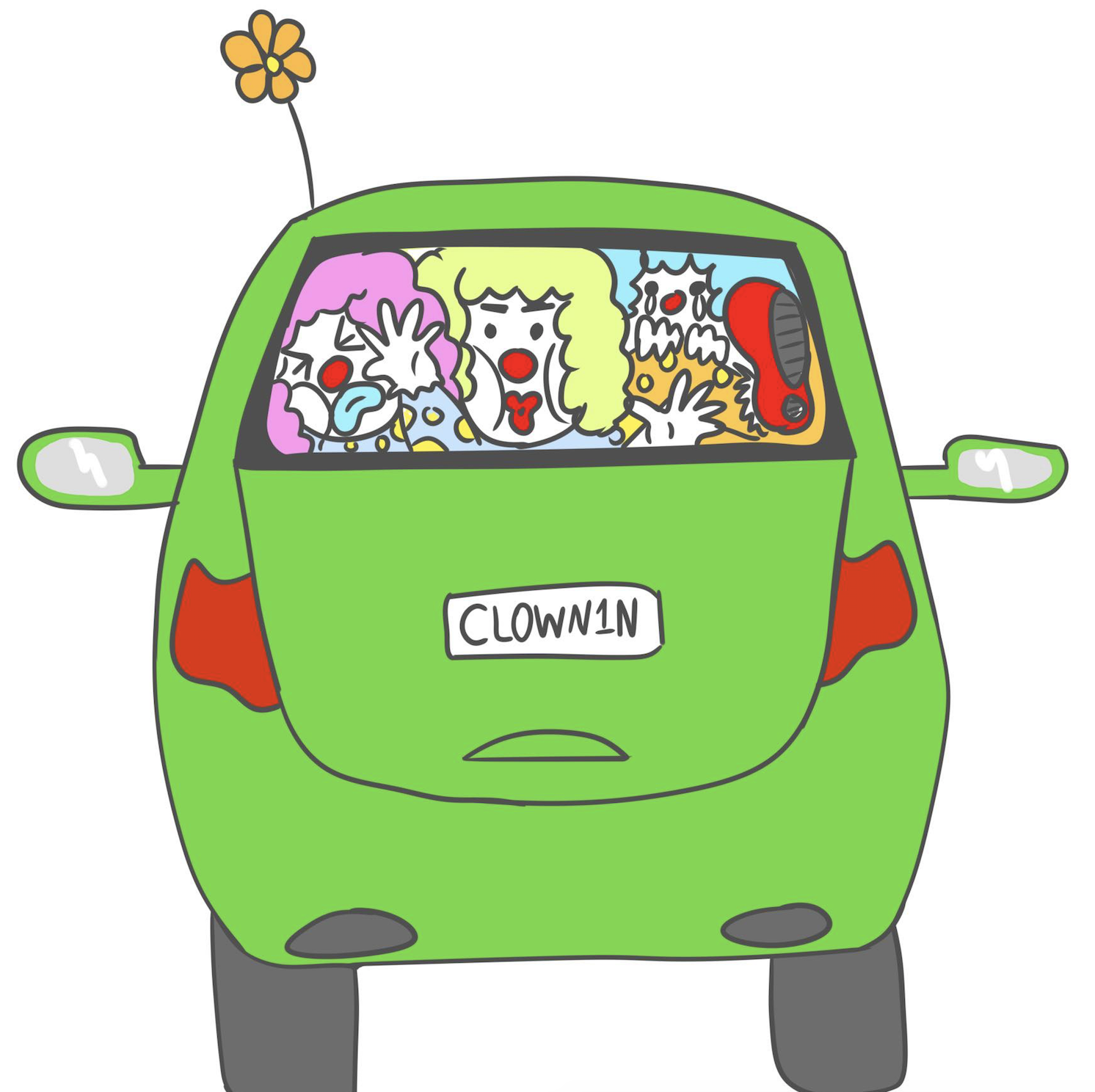}
\end{wrapfigure} 

 This problem requires you to put on your thinking cap and your clown shoes. Imagine there are $n\in\N$ clowns that need to go to the circus. Suppose they line up to board $m\in\N$ clown cars that will take them to the circus tent. A number $d \in \N$ indicates the number of clowns that fit in a car. This number $d$ can be fixed or varied, whichever you decide. Now, suppose each clown has a preferred clown car, which is represented in a preference list. If a clown gets to their preferred clown car and it is full, i.e. there are already $d$ clowns in it, the clown will continue forward, checking clown cars until it finds one with a spot available. If all the clowns fit in the $m$ cars, they are all able to make it to the circus. If all the clowns make it to the circus based on the given preference list, the list of clown preferences is a clown function. For given integers $d$, $m$, and $n$, can you count all clown functions?

Finished? Good news! This isn't the end of the road! Once you have counted these parking functions, there are many more properties that you can study. Go to \bref{hungryformore}. 

\subsection{--\ \ Restricted backward movement}\label{fixeddirection}
$ $\newline

Even the most expensive cars require maintenance. Suppose that out of the $n$ cars, $c_1, c_2, \dots, c_n$, every car with an even index has a dysfunctional gear shift preventing it from being able to reverse. Thus, cars $c_2, c_4, c_6$, etc., can only move forward to check the spots ahead of them if their preferred spot is occupied. However, cars at odd indices, $c_1, c_3, c_5$, etc., have no outstanding issues and are equipped with fully functioning gear shifts. This means that, if their preferred spot is occupied, they first check one spot behind them and park there if it is empty. Otherwise, they continue forward and check the remaining spots in the lot, parking in the first available. In this scenario, find the number of preference lists that allow for all of the cars to park.

Finished? Good news! This isn't the end of the road! Once you have counted these parking functions, there are many more properties that you can study. Go to \bref{hungryformore}. 

\subsection{--\ \ Introducing randomness}\label{randomdirection}
$ $\newline
\begin{wrapfigure}{r}{1.7in}
\begin{tikzpicture}
\node at (-.25,.25){\includegraphics[width=1.5in]{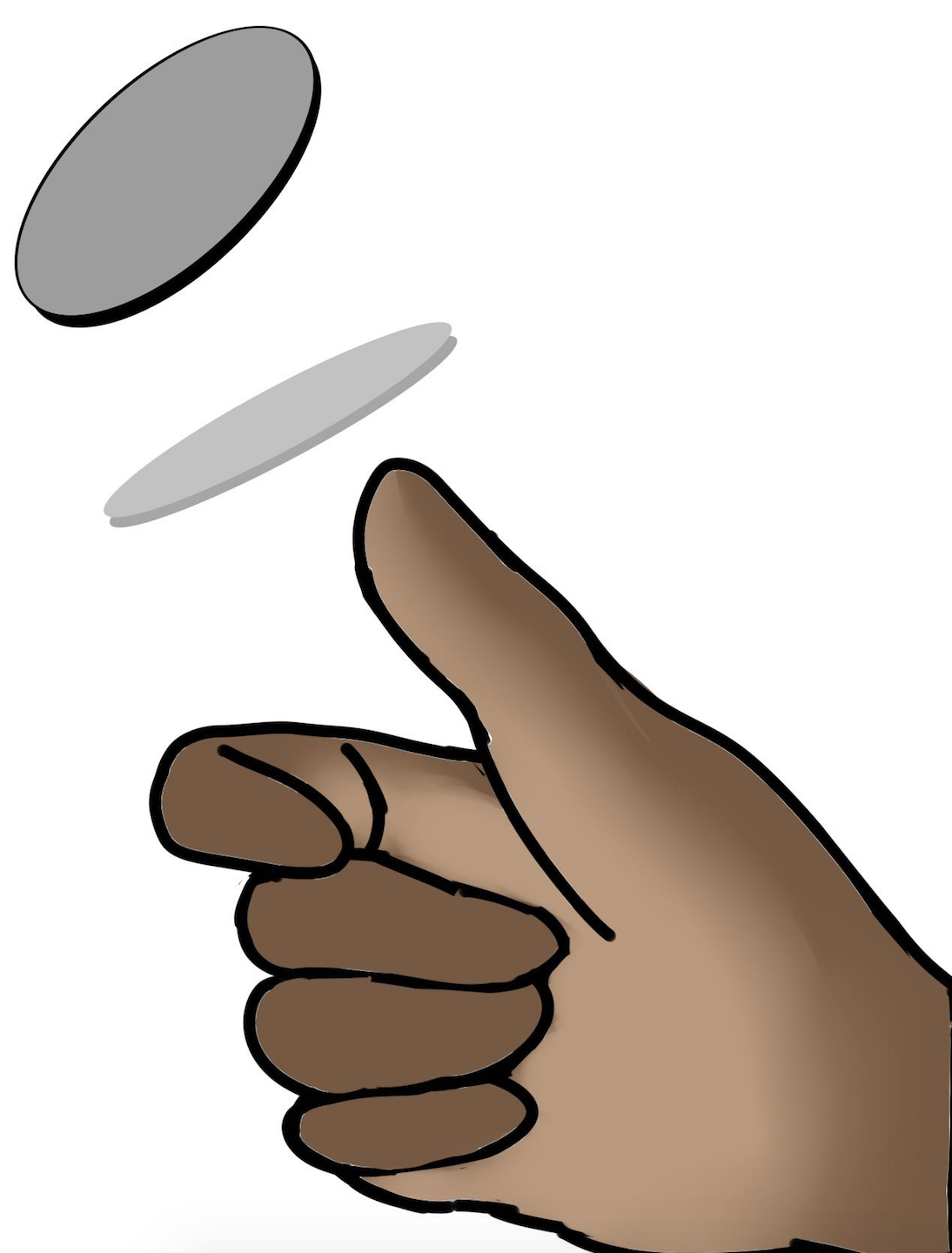}};
\node at (0,2.5){\large{Heads}};
\node at (.5,1.75){\large{or}};
\node at (1,1.0){\large{tails?}};
\end{tikzpicture}
\end{wrapfigure} 

Making decisions can be hard sometimes. The drivers in this parking scenario agree. Suppose that upon discovering their preferred spot is taken, the driver of each car flips a fair coin to decide whether they will back up and check the previous spot or if they will just continue forward to find the next available spot. If a car checks the previous spot and finds that it is taken, they continue driving forward to look for the next unoccupied spot. Given a fixed preference list, what is the probability that the preference list is a parking function?  

Finished? Good news! This isn't the end of the road! Once you understand these probabilities, there are many more properties that you can study. Go to \bref{hungryformore}.

\subsection{--\ \ Obstructed parking spots}\label{backwardtrailer}$ $ \newline

Imagine a neighborhood where every family decides to host a BBQ for their friends on the same day. Suppose that there are $n \in \N$ cars arriving and they all need to park. The neighborhood street has $m \in \N$ spots with $m > n$.
Some of the $m$ parking spots are obstructed by grills, lawn games, and other BBQ festivities for the duration of the event, limiting the number of available spots in which the BBQ attendees can park. Some hosts are more popular than others and not every BBQ is the same size. This means that some hosts take up more spots in the parking lot along the neighborhood street than others. Luckily, there are still $n$ spots available, one for every car who wants to park in the lot. However, they must park around these obstructions.

Suppose that each car $c_i$ has a preferred parking spot, $a_i$, with $1\leq a_i \leq m$ that is closer to the BBQ they want to attend. If their preferred spot is empty when they arrive, they park in that spot. However, their preferred spot could be taken by another guest that arrived earlier or be obstructed by a BBQ activity. If the car's preferred spot is occupied, the car can back up one spot. If it is empty, the relieved guest will take it. If it is occupied, the car will continue driving forward until it finds an empty spot. If no such spot exists, much to the dismay of the excited guest, they will return home. The preference list consisting of the preferences of these cars is a parking function if all the guests can park and attend their BBQ. How many such parking functions are there?

Finished? Good news! This isn't the end of the road! Once you have counted these parking functions, there are many more properties that you can study. Go to \bref{hungryformore}.

\subsection{--\ \ Cars and scooters}\label{backwardmotorcyle} 

$ $\newline

\begin{wrapfigure}{l}{2in}
\includegraphics[trim=2mm 0 0 0, clip,width=2in]{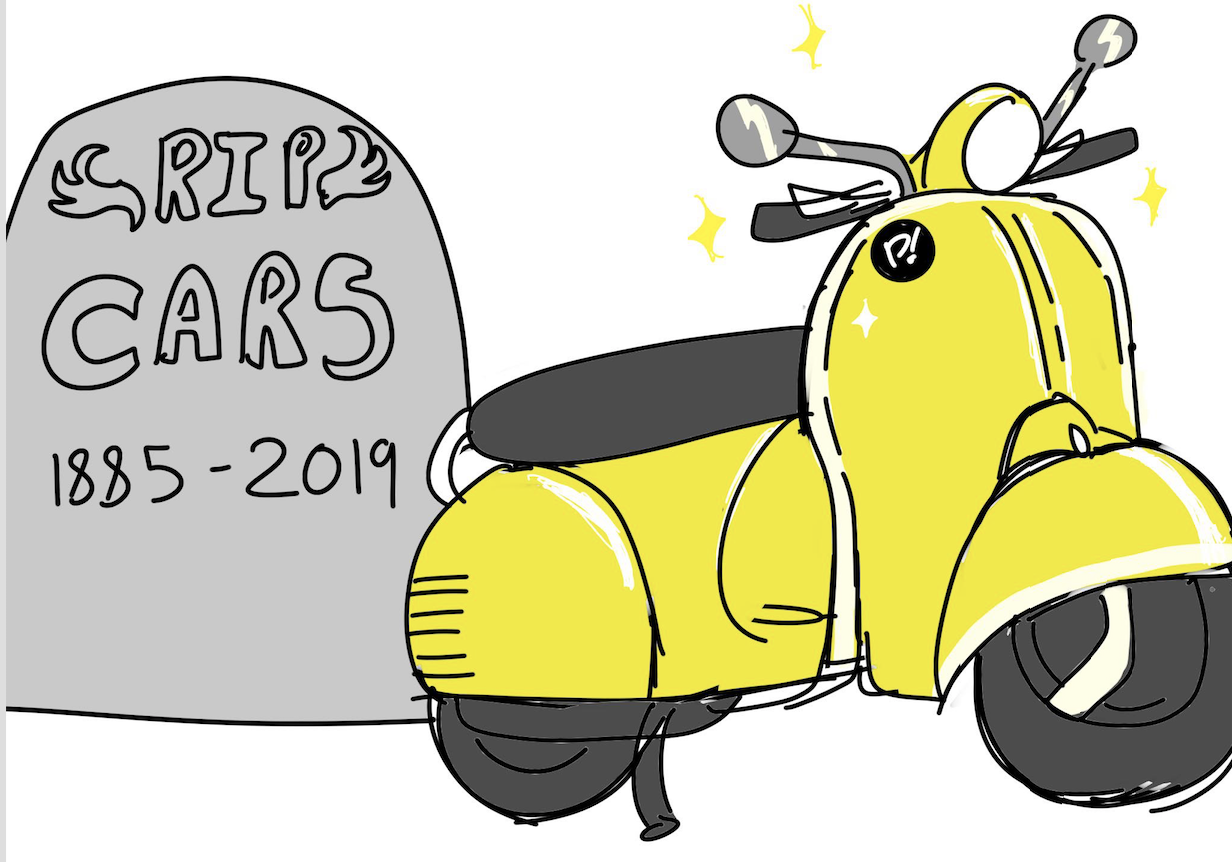}\\
\end{wrapfigure} 

Imagine that you live in a city where everyone is trading in their cars for sleeker, more space efficient vehicles: electric scooters. Suppose that there are $n \in \mathbb{N}$ scooters lined up to park in $m \in \mathbb{N}$ parking spots. Each scooter has a preferred parking spot. Since scooters are smaller than cars, many scooters can park in a single parking spot and $m\leq n$. We denote the number of scooters that can fit in a parking spot by $d$. You get to determine the value of $d$ however you see fit. If a scooter gets to its preferred parking spot and it is already full of scooters, it can check the previous spot. If the previous spot has room, then the scooter will park there. If not, the scooter will continue forward until it finds a spot that has not yet reached max capacity, in which case it parks there. Can you count the number of scooter parking functions for an arbitrary $d$, $m$, and $n$?

Finished? Good news, this is not the end of the road! Once you have counted these parking functions, there are many more properties that you can study. Go to \bref{hungryformore}.

\subsection{--\ \ Teleporting cars}\label{teleportation} $ $ \newline
 \begin{wrapfigure}{r}{2in}
\includegraphics[width=2in]{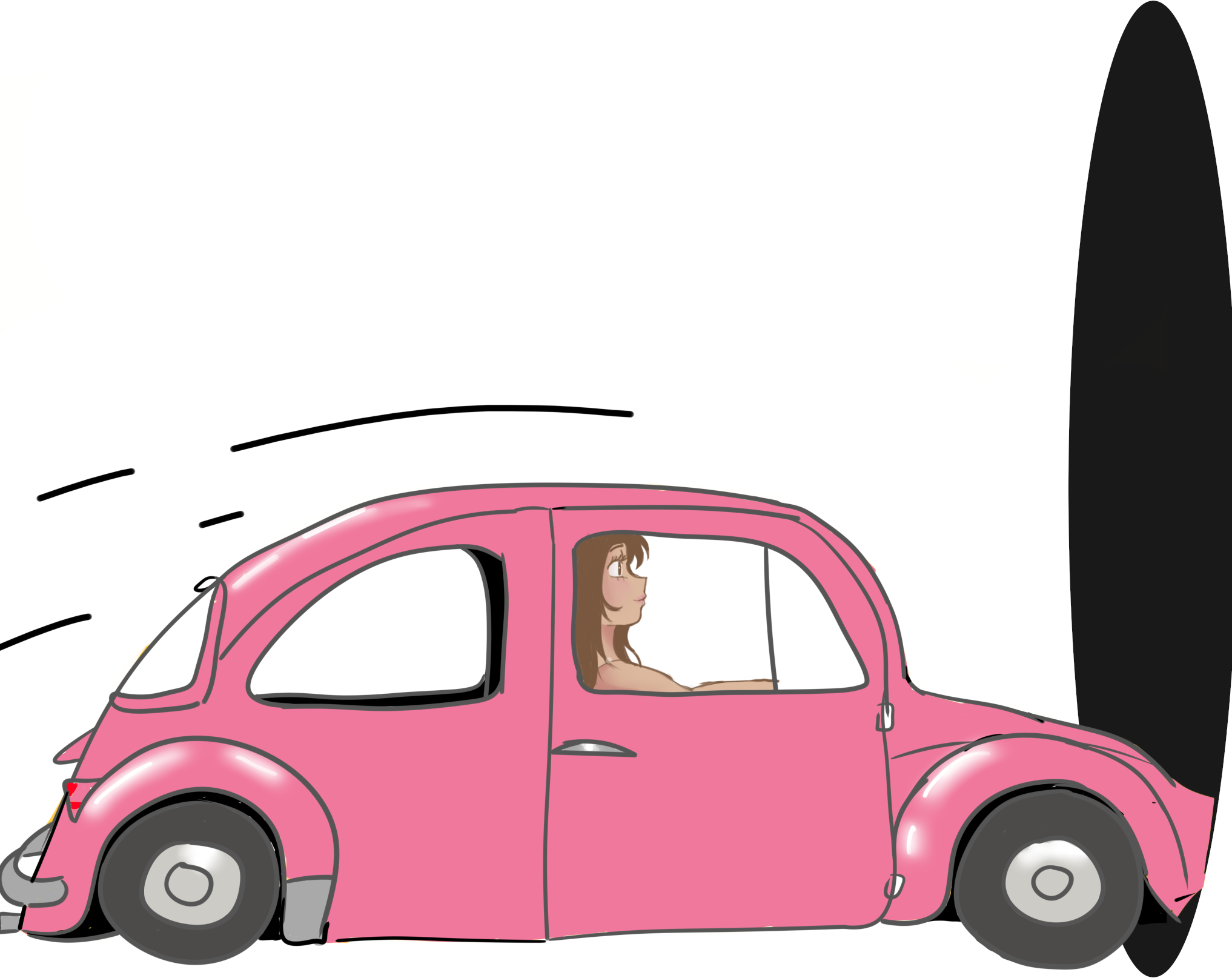}
\end{wrapfigure} 
 Oh no! Something wonky is going on and the parking lot is now riddled with wormholes. As usual, there are $n$ cars that want to park in $n$ parking spots. However, in this scenario, if car $c_i$'s preferred parking spot, $a_i$, is occupied, it falls into a wormhole that teleports it $k$ parking spots backward. You can choose how the value of $k$ is determined. After teleporting, the car remains intact and decides to take advantage of the worm hole situation by checking the spot, $a_i -k$, that it landed next to. If that spot is empty, the car will park there. If not, it proceeds forward as usual until it finds an empty parking spot.
 
 If all cars are able to park based on a given preference list, then the preference list is a parking function. If a car teleports outside of the $n$ parking spots or is unable to park, then the preference list is not a parking function. For a given $k,n\in\N$, how many preference lists are parking functions?

 Finished? Good news! This isn't the end of the road! Once you have counted these parking functions, there are many more properties that you can study. Go to \bref{hungryformore}.

 \subsection{--\ \ Back of the line blues}\label{i-stepsback} $ $\newline 

In an ideal world, every car would have a fair chance at finding a spot in the parking lot. Unfortunately, life is not always fair and in this problem, the cars at the end of the preference list are at a disadvantage. Suppose the number of steps back a car can take is determined by its position in line. In this scenario, cars at the beginning of the line can take more steps back than cars at the end. Specifically, define $f: [n] \rightarrow [n]$ as $f(i) = n-i+1$. For each $i \in [n]$, if car $c_i$'s preferred spot $a_i$ is taken, it can check up to $f(i)$ spots behind $a_i$ before continuing forward. Note that for each pair of indices $i,j$ with $i<j$, $f(j) < f(i)$. With this rule in place, how many preference lists allow all the cars to park?

Finished? Good news! This isn't the end of the road! Once you have counted these parking functions, there are many more properties that you can study. Go to \bref{hungryformore}. 

\section{The adventure continues}\label{hungryformore} 

Now that you have finished enumerating your parking problem, how do you feel? We will let these feelings guide the next part of your adventure. 
Are you feeling... 
\begin{itemize}
    \item[a.] Uncertain? Go to \bref{IntroducingRandomness}.
    \item[b.] Sad to move on? See \bref{Bumping}.
    \item[c.] On top of the world? Go to \bref{PeaksandValleys}.
    \item[d.]  Like it is all downhill from here? Go to \bref{AscentsTiesDescents}.
    \item[e.] Average? Go to \bref{Distributions}.
    \item[f.] Lucky? Go to \bref{LuckyCars}.
\end{itemize}

\subsection{--\ \ Embracing uncertainty
}\label{IntroducingRandomness} $ $ \newline

If you haven't already considered adding an aspect of randomness to your problem, see the problem presented in \bref{randomdirection} for an idea of how to do so. If you have, consider a weighted coin or randomizing a different or additional decision in the parking process.

These suggestions alone can set you up for a lifetime of adventure. However, if you are hoping to go in another direction, return to \bref{firstQuestion} and see where a new path can take you. When you are ready to conclude your journey, go to \bref{Conclusion}.

\subsection{--\ \ High fives and fist bumps}\label{Bumping} $ $ \newline 

 \begin{wrapfigure}{l}{2in}
\includegraphics[trim=0 70mm 0 0, clip,width=2in]{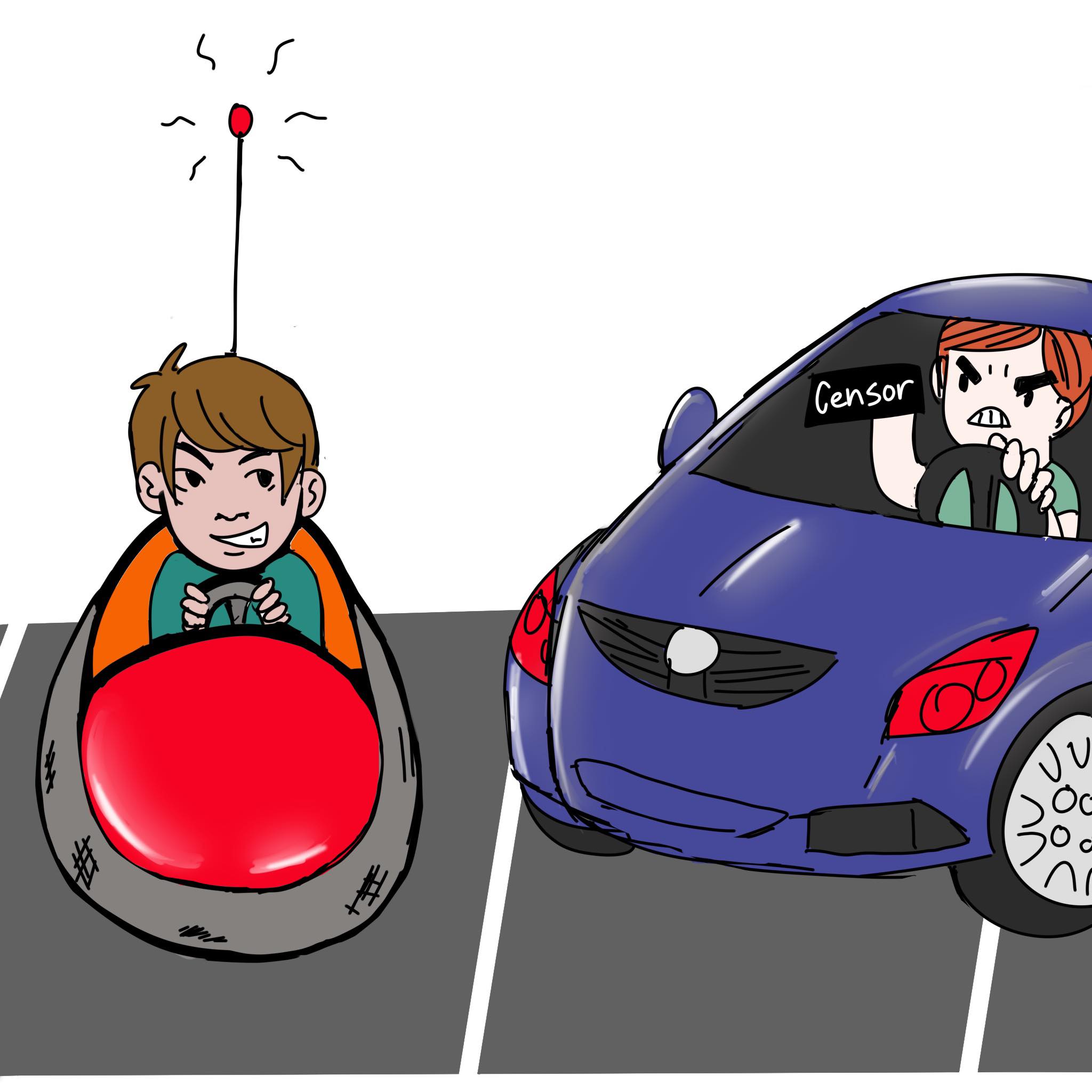}
\end{wrapfigure}

Moving on from a good problem can be challenging. A car moving on from its preferred parking spot can be challenging as well. As the cars park, they do not always get their preferred spot since it may be occupied. In this case, they too feel sad to move on. This sadness is something you can study. When a car arrives at a parking spot and is forced to move on, i.e. check another spot because the current one is full, it is called a bump. For any parking function, the number of bumps can be counted as the distance between where a car wants to park and where it has to park, which we call displacement. Studying the displacement of these parking functions leads to a new research project.

These suggestions alone can set you up for a lifetime of adventure. However, if you are hoping to go in another direction, return to \bref{firstQuestion} and see where a new path can take you. When you are ready to conclude your journey, go to \bref{Conclusion}.

\subsection{--\ \ The highs and lows of parking functions}\label{PeaksandValleys} $ $ \newline

So you're feeling on top of the world. That's awesome, you should be feeling that way! Some elements of a preference list feel on top of their little world too. When a number in a preference list feels like this, it is called a peak. Peaks occur at an index $i$ when car $c_i$ has a preference greater than the preferences of the two cars ($c_{i-1}$ and $c_{i+1}$) next to it. Specifically, in the preference list $\alpha=(a_1,a_2,\ldots,a_n)$, a peak occurs when $a_{i-1} <a_i$ and $a_{i+1} < a_i$ for some $i \in \{2, 3, \ldots, n-1\}$. How many peaks are there in your particular parking functions? What is the total number of peaks in the parking functions with $n$ cars? These are the types of questions you can ask to continue on in your adventure. 

However, if you get stuck in a rut while thinking about peaks, don't be discouraged. Valley's are another property you can try to count in your scenario. A valley is the opposite of a peak and occurs when car $c_i$ has a preference less than the preferences of the two cars ($c_{i-1}$ and $c_{i+1}$) next to it. 

These suggestions alone can set you up for a lifetime of adventure. However, if you are hoping to go in another direction, return to \bref{firstQuestion} and see where a new path can take you. When you are ready to conclude your journey, go to \bref{Conclusion}.

\subsection{--\ \ What goes up must come down, unless it's all tied up}\label{AscentsTiesDescents} $ $ \newline

Although you feel like it is all downhill from here, you can always turn a corner and take a path that leads you back up or keeps you leveled. Luckily, no matter where your journey takes you, your parking function problem involves preference lists and you can explore and study their characteristics. 

To describe these characteristics, we consider the preference list $\alpha=(a_1,a_2,\ldots,a_n)$, and for all $1\leq i\leq n-1$, we~say $i$ is an ascent if $a_i<a_{i+1}$, 
    $i$ is a descent if $a_i>a_{i+1}$,  and
    $i$ is a tie if $a_i=a_{i+1}$.
You can count these for any parking function, regardless of the parking rules.

Studying these characteristics in your parking problem can set you up for a lifetime of adventure. However, if you are hoping to go in another direction, return to \bref{firstQuestion} and see where a new path can take you. When you are ready to conclude your journey, go to \bref{Conclusion}.

\subsection{--\ \ Probabilistic-alley}\label{Distributions} $ $ \newline

If you're feeling average, get ready to put on your probabilistic hat. Given a random parking function, what is the probability that a given car will get its preferred parking spot? Among the total number of ascents and descents\footnote{For the definition of ascents and descents, go to \bref{AscentsTiesDescents}, but don't get too distracted.} in the set of all parking functions, how are they distributed? That is, are certain parking spots more likely than others to have ascents or descents in the preference list?

Finding the answers to these questions can set you up for a lifetime of adventure. However, if you are hoping to go in another direction, return to \bref{firstQuestion} and see where a new path can take you. When you are ready to conclude your journey, go to \bref{Conclusion}.

\subsection{--\ \ 
Count your lucky stars
}\label{LuckyCars} $ $ \newline
\begin{wrapfigure}{r}{2in}
\includegraphics[trim=0 0 0 0, clip,width=2in]{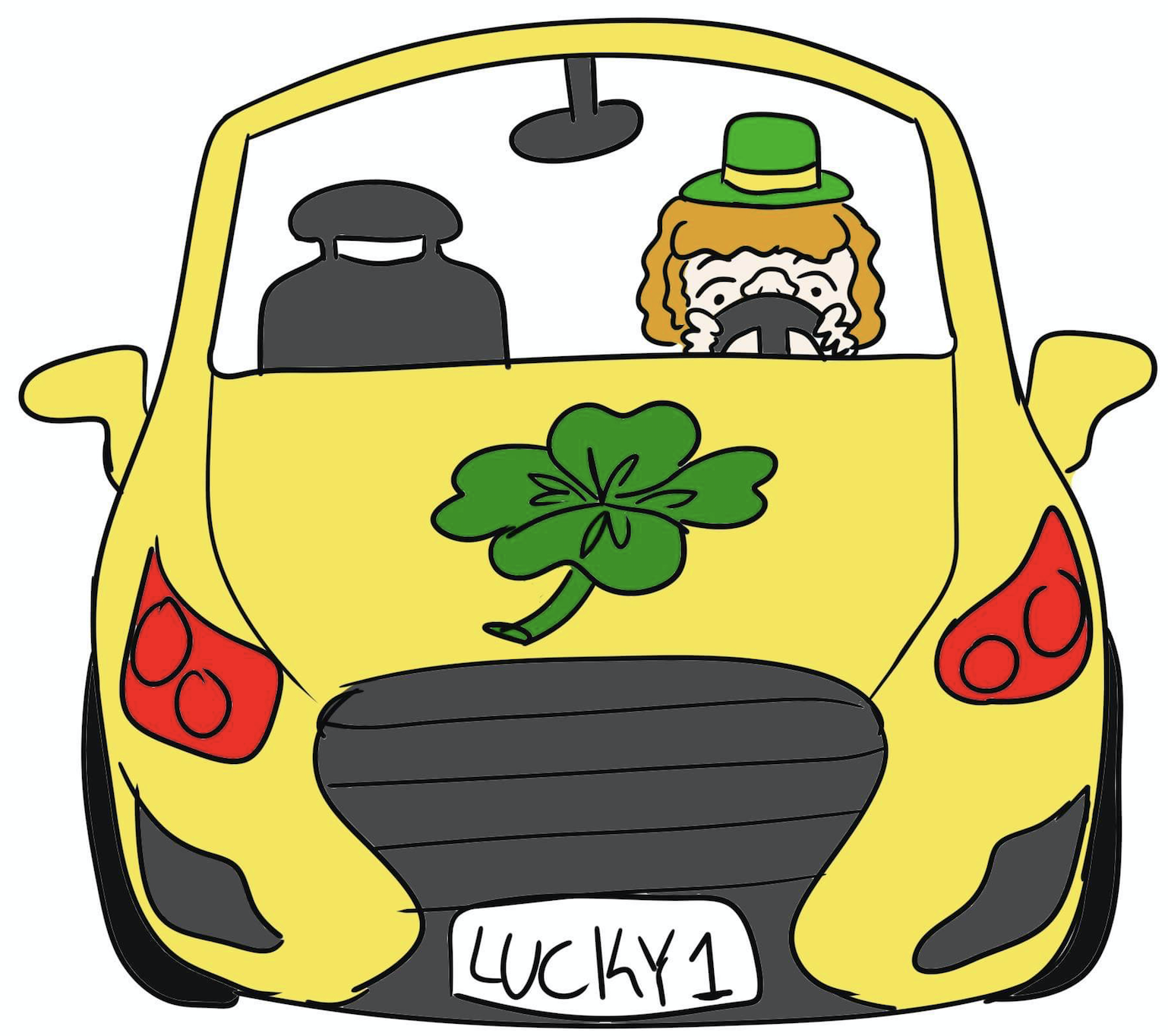}
\end{wrapfigure} 
All $n$ cars are able to park when the preference list $(a_1, a_2, \ldots, a_n)$ is a parking function. However, some cars may not get to park exactly in their preferred spot. Car $c_i$ is defined to be lucky if it parks in spot $a_i$. Try counting the number of lucky cars for your parking problem.

This suggestion alone can set you up for a lifetime of adventure. However, if you are hoping to go in another direction, return to \bref{firstQuestion} and see where a new path can take you. When you are ready to conclude your journey, go to \bref{Conclusion}.

\section{The end of the beginning}\label{Conclusion}
Now let's pump the brakes for a second. As you can imagine, there are many other combinations of these different parking function scenarios. With so many variations, if you reached the end of this paper and did not find a parking function problem that tooted your horn, we encourage you to build your own.

\bigskip

\begin{center}
    {\huge\calligra{Your adventure awaits you!}}
\end{center}

\section*{Epilogue}\label{sec:closingremarks}
In this short adventure, we simulate for the reader the twists and turns that encompass the research process. Sometimes, we cruise through obstacles with maximum productivity. Other times, we seem to exhaust our resources and it feels like no more progress can be made. Ultimately, the journey is never over as we continue to find new paths to explore. These paths can sometimes take us down unexpected roads, which can lead to interesting connections. Many researchers have experienced this. Thus, the connections of parking functions to other fields are plenty. For instance, parking functions have been studied on graphs in \cite{treeParking}. Alternatively, in \cite{PitmanStanley}, parking functions have been connected to volumes of Pitman-Stanley polytopes. Moreover, generalizations abound and all that is left is more time to enumerate them and discover their properties. We are counting on you to continue on the quest to solve parking problems!

\begin{table}[h]
    \centering
    \resizebox{.94\textwidth}{!}{
    \begin{tabular}{|>{\arraybackslash}m{2.25cm}|m{5.55in}|}
         \multicolumn{2}{c}{{\bf{These adventures were made possible in part by$\ldots$}}}
         \\[5pt]\hline
         \bref{onlyforward}, \bref{carsizes}& Cars with variable size have been studied in \cite{ehrenborg1} and allowing multiple cars to park in a spot have been studied in \cite{moto}. However they did not include variable movement, nor any other generalization presented here.
         \\\hline
         \bref{backwardthenforward}, \bref{onespaceback}, 
         \bref{multiplespacesback}
         &  Cars that can back up to a certain number of spots when finding their preferred parking spot occupied were studied in \cite{Ours}. This only included looking one parking spot back at a time. It did not include moving back a fixed number of spots and taking the first available.\\\hline
        
         \bref{variablemovement}& Cars whose backward movement varies depending on the car is an open problem. \\\hline
         
         \bref{onespacemanycars}, \bref{fixeddirection}, \bref{randomdirection}, \bref{teleportation},
         \bref{i-stepsback},
         \bref{IntroducingRandomness}& The problems in these sections have not been studied. \\\hline
         \bref{backwardtrailer}& 
Cars parking after a trailer were studied in \cite{ehrenborg2}. More recently, in \cite{preprintPH}, the authors generalize these results to cars parking after a fixed number of trailers have already been parked, and where cars have only forward movement. Thus cars with variable or backward movement have not been studied in these scenarios.
\\\hline
         \bref{backwardmotorcyle}& A variation of this problem has been studied in \cite{moto} where the authors make motorcycles indistinguishable and cars distinguishable. All vehicles only moved forward. Allowing for variable movement leads to open problems.\\\hline
         \bref{Bumping}& For the study of displacement for classical parking function see \cite{kreweras,yan}.  \\\hline
                 \bref{PeaksandValleys}& This problem has been studied on permutations \cite{BBS}. An initial study for these results for parking functions that are not permutation was done in \cite{MSM}. No further study has been made on generalized parking functions.
\\\hline
         \bref{AscentsTiesDescents}& The enumeration of descents and ties of classical parking functions was studied in \cite{pfdt}. For more examples of what has been studied about ascents, descents, and ties in parking functions, see \cite{TextbookPam,ProbabilizingParkingFunctions,pfdt}. \\\hline
 \bref{Distributions}& Distributions and probabilistic problems on classical parking functions were studied in \cite{ProbabilizingParkingFunctions}. This has not been expanded to other generalized parking functions.
 \\\hline
         \bref{LuckyCars}& Lucky cars were defined and studied in \cite{yan} and a generating function for the number of lucky cars in a parking function is obtained in \cite{GesselSeo}. This has not been expanded to other generalized parking functions. \\\hline
    \end{tabular}
    }
\end{table}

\noindent{\bf{Acknowledgements}}.
This research was supported in part by the Alfred P. Sloan Foundation, the Mathematical Sciences Research Institute, and the National Science Foundation. We thank Ana Valle for all of the artwork. 

\addresseshere

\bibliographystyle{numeric}
\bibliography{ref}

\end{document}